\documentclass{amsart}

\usepackage{amsfonts,amsmath,amssymb,amsthm,mathrsfs,amsxtra}
\usepackage{enumerate,verbatim}
\usepackage{lineno} 
\usepackage[all,2cell,ps]{xy}
\usepackage[pagebackref]{hyperref}

\DeclareMathOperator{\ann}{ann}

\DeclareMathOperator{\curv}{curv}
\DeclareMathOperator{\cx}{cx}
\DeclareMathOperator{\depth}{depth}

\DeclareMathOperator{\Ext}{Ext}

\DeclareMathOperator{\gdim}{G-dim}

\DeclareMathOperator{\Hom}{Hom}

\DeclareMathOperator{\injcx}{inj\,cx}
\DeclareMathOperator{\injdim}{injdim}
\DeclareMathOperator{\pd}{pd}
\DeclareMathOperator{\projdim}{pd}

\DeclareMathOperator{\px}{px}
\DeclareMathOperator{\rank}{rank}
\DeclareMathOperator{\Soc}{Soc}
\DeclareMathOperator{\Tor}{Tor}

\DeclareMathOperator{\type}{type}

\renewcommand{\ge}{\geqslant}
\renewcommand{\le}{\leqslant}

\newcommand{\fm}{\mathfrak{m}}

\renewcommand{\iff}{if and only if }

\newcommand{\oM}{\overline{M}}
\newcommand{\oN}{\overline{N}}
\newcommand{\oomega}{\overline{\omega}}
\newcommand{\oR}{\overline{R}}

\theoremstyle{plain}
\newtheorem{theorem}{Theorem}[section]
\newtheorem{lemma}[theorem]{Lemma}
\newtheorem{proposition}[theorem]{Proposition}
\newtheorem{corollary}[theorem]{Corollary}

\theoremstyle{definition}
\newtheorem{definition}[theorem]{Definition}
\newtheorem{example}[theorem]{Example}

\newtheorem{question}[theorem]{Question}

\newtheorem{para}[theorem]{}

\theoremstyle{remark}
\newtheorem{remark}[theorem]{Remark}

\numberwithin{equation}{theorem}

\title[Gorensteinness of short local rings]{Gorensteinness of short local rings in terms of the vanishing of Ext and Tor}
\author[Dipankar Ghosh]{Dipankar Ghosh}
\address{Department of Mathematics, Indian Institute of Technology Hyderabad, Sangareddy, Kandi, Telangana - 502285, India}
\email{dghosh@math.iith.ac.in, dipug23@gmail.com}
\date{August 25, 2020}

\begin{document}


\begin{abstract}
	Let $(R,\fm)$ be a commutative Noetherian local ring which contains a regular sequence $ \underline{x} = x_1,\ldots,x_d \in \fm \smallsetminus \fm^2 $ such that $ \mathfrak{m}^3 \subseteq (\underline{x}) $. Let $ M $ be a finite $ R $-module with maximal complexity or curvature, e.g., $ M $ can be a nonzero direct summand of some syzygy module of the residue field $ R/\fm $. It is shown that the following are equivalent: (1) $R$ is Gorenstein, (2) $\Ext_R^{\gg 0}(M,R)=0$, and (3) $\Tor_{\gg 0}^R(M,\omega) = 0$, where $\omega$ denotes a canonical module of $R$. It gives a partial answer to a question raised by Takahashi. Moreover, the vanishing of $\Ext_R^{\gg 0}(\omega,N)$ for certain $ R $-module $ N $ is also analyzed. Finally, it is studied why Gorensteinness of such local rings is important.
\end{abstract}
\maketitle

\section{Introduction}\label{sec:intro}
 
 Let $ R $ be a local ring with the maximal ideal $ \fm $ and residue field $ k $. For every $ n \ge 0 $, $ \Omega_n^R(k) $ denotes the $ n $th syzygy module in a minimal free resolution of $ k $.
 
 One of the most influential results in commutative algebra is the result of Auslander, Buchsbaum and Serre:  $ R $ is regular if and only if projective dimension of $ k $ is finite, which is equivalent to the fact that some syzygy module of $k$ is free. Dutta, in \cite[Cor.~1.3]{Dut89}, proved that $R$ is regular if and only if $\Omega_n^R(k)$ has a nonzero free direct summand for some $ n \ge 0 $. Later, in \cite[Thm.~4.3]{Tak06}, Takahashi generalized Dutta's result by showing that $ R $ is regular if and only if $ \Omega_n^R(k) $ has a semidualizing (e.g., $ R $ itself as an $ R $-module) direct summand for some $ n \ge 0 $. In a joint work \cite[Cor.~3.2]{GGP} with Gupta and Puthenpurakal, the author proved a considerably stronger result: If a finite direct sum of syzygy modules of $k$ maps onto a semidualizing $R$-module, then $R$ is regular.
 
 It follows from Dutta's result that $ R $ is regular if and only if some $ \Omega_n^R(k) $ ($ n \ge 0 $) has a nonzero direct summand of finite projective dimension. A counterpart of this result for injective dimension is shown in \cite[Thm.~3.7]{GGP}. In a different direction, Martsinkovsky \cite[Prop.~7]{Mar96} generalized Dutta's result as follows: If a finite direct sum of syzygy modules of $ k $ maps onto a nonzero $ R $-module of finite projective dimension, then $ R $ is regular. Thereafter, Avramov proved a much more stronger result:
 
 \begin{theorem}[Avramov]{\rm \cite[Cor.~9]{Avr96}}\label{thm:Avr-max-cx}
 	Every nonzero homomorphic image $ M $ of a finite direct sum of syzygy modules of $ k $ has maximal complexity and curvature, i.e., $ \cx_R(M) = \cx_R(k) $ and $ \curv_R(M) = \curv_R(k) $; see {\rm Definition~\ref{def:complexity}}.
 \end{theorem}

\begin{para}\label{para:proj-CI-dim}
	A few consequences of Theorem~\ref{thm:Avr-max-cx} are the following results. Let $ M $ be a nonzero homomorphic image of a finite direct sum of syzygy modules of $ k $. (1) If projective dimension $ \projdim_R(M) $ is finite, then $ R $ is regular; one may use Lemma~\ref{lem:cx-properties}(i)(a). (2) If complete intersection dimension (as in \cite[(1.2)]{AGP97}) of $ M $ is finite, then $ R $ is a complete intersection ring; use \cite[(5.6)]{AGP97} and Proposition~\ref{prop:facts-projcx-k}(ii).
\end{para}

The Gorenstein dimension (in short G-dimension) was introduced by Auslander \cite{Aus67}, and was deeply studied by him and Bridger \cite{AB69}. It is well known that $\gdim_R(k)$ is finite if and only if $R$ is Gorenstein. In this theme, Takahashi showed that $R$ is Gorenstein if and only if $ \Omega_n^R(k) $ has a nonzero direct summand of G-dimension $ 0 $ for some $ 0 \le n \le \depth(R) + 2 $; see \cite[Thm.~6.5]{Tak06}. The following question \cite[6.6]{Tak06} of Takahashi is still open: If $ \Omega_n^R(k) $ has a nonzero direct summand of G-dimension $ 0 $ for some $ n > \depth(R) + 2 $, then is $R$ Gorenstein? More generally, in view of Section~\ref{para:proj-CI-dim}, it is now natural to ask the following:

\begin{question}\label{ques1}
	If a finite direct sum of syzygy modules of $ k $ maps onto a nonzero $ R $-module of finite G-dimension, then is $R$ Gorenstein?
\end{question}

We give a positive answer to Question~\ref{ques1} for short local rings:

\begin{definition}\label{short-local-ring}
	The ring $ R $ is said to be a {\it short local ring} if it contains a regular sequence $ \underline{x} := x_1,\ldots,x_d \in \fm \smallsetminus \fm^2 $ with the property that $ \mathfrak{m}^3 \subseteq (\underline{x}) $.
\end{definition} 

The motivation to work over these rings came from the following results. A commutative local analog of a conjecture of Tachikawa says that if $ \Ext_R^i(\omega, R) = 0 $ for all $ i \ge 1 $, then $ R $ is Gorenstein, where $ \omega $ is a canonical $ R $-module. In \cite{ABS05}, Avramov, Buchweitz and \c{S}ega proved this conjecture in several significant cases. In a particular case, they considered short local rings;
see \cite[Thm.~5.1]{ABS05}. Many important conjectures that are still open in general have been verified over these rings; see \cite[Thm.~B]{Les85} and \cite[Thm.~2.11 and 4.1]{HSV04}. As it is described in \cite[page~459]{AIS08}, there are many famous counterexamples too built over these rings. In the present article, inspired by Question~\ref{ques1}, we prove the following:
 
 \begin{theorem}[=~\ref{cor:Gor-rings-hom-image-of-syz}]\label{thm1}
 	Let $ (R,\mathfrak{m},k) $ be a short local ring. Let $ \omega $ be a canonical module of $ R $.
 	Let $ M $ be a nonzero homomorphic image of a finite direct sum of syzygy modules of $ k $. Then the following statements are equivalent:
 	\begin{enumerate}[\rm (i)]
 		\item $ R $ is Gorenstein.
 		\item $ \Ext_R^i(M,R) = 0 $ for all $ i \gg 0 $.
 		\item $ \Tor_i^R(\omega, M) = 0 $ for all $ i \gg 0 $.
 	\end{enumerate}
 \end{theorem}
 
 With hypotheses as in Theorem~\ref{thm1}, we also study the following:
 
 \begin{question}\label{ques2}
 	If $ \Ext^i_R(\omega, M) = 0 $ for all $ i \gg 0 $, then is $ R $ Gorenstein?
 \end{question}

 Being $ \gdim_R(M) $ finite is much stronger condition than $ \Ext_R^{\gg 0}(M,R) = 0 $;
 see \cite[Remarks after (3.7)]{AB69}. Hence Theorem~\ref{thm1} ensures that Question~\ref{ques1} has affirmative answer for short local rings.
 If $ R $ has minimal multiplicity (i.e., if $ \fm^2 = (\underline{x}) \fm $ for some $ R $-regular sequence $ \underline{x} $) and $ M $ ($ \neq 0 $) is a direct summand of some syzygy module of $ k $, then Question~\ref{ques2} has positive answer; see \cite[Thm.~5.9]{GP} and \cite[Thm.~5.5]{Gho}. In this article, though we are unable to give a complete answer to Question~\ref{ques2}, but in various attempts, we prove Theorems~\ref{thm:Gor-rings-max-injcx}, \ref{cor:length-neq-2e-Gor-rings} and \ref{prop:mu-le-length-mM-vanishing-of-Ext}.
 
 We now describe in brief the contents of this article. In Section~\ref{sec:preliminaries}, we collect some preliminaries for later use. Our main results are shown in Section~\ref{sec:main}, which provide necessary and sufficient conditions for a short local ring to be Gorenstein in terms of the vanishing of certain Ext and Tor modules. Finally, in Section~\ref{sec:last}, we give some reason why Gorensteinness of such local rings is important; see Theorem~\ref{thm:char-of-reg} and Example~\ref{exam1:thm:char-of-reg}.
 
\section{Preliminaries}\label{sec:preliminaries}
 
 Throughout this article,
 all rings are assumed to be commutative Noetherian local rings, and all modules (except possibly the injective hulls) are assumed to be finitely generated. Moreover, {\it $ R $ always denotes a local ring with the unique maximal ideal $\mathfrak{m}$ and residue field $k$}. For an $ R $-module $ M $, and integer $ n \ge 0 $, $\Omega_n^R(M)$ denotes the $ n $th {\it syzygy module} of $ M $ in a minimal free resolution. Though $\Omega_n^R(M)$ depends on the choice of a minimal free resolution of $ M $, but it is unique up to isomorphism. For $ n \ge 1 $, since $ \Omega_n^R(M) \subseteq \mathfrak{m} F $ for some free module $ F $, one obtains the following relation between the socle of the ring and the annihilator of the syzygy modules:
 
 \begin{lemma}\label{lemma: socle, syzygy}
 	$ \Soc(R) \subseteq \ann_R\left(\Omega_n^R(M)\right) $ for every $ n \ge 1 $.
 \end{lemma}
 
 We set $ M^* := \Hom_R(M,R) $ and $ M^{\vee} := \Hom_R(M,E) $, where $ E := E_R(k) $ is the injective hull of $ k $ over $ R $. We denote the minimal number of generators of $ M $ by $ \mu(M) $, and the length of $ M $ by $ \lambda(M) $. For every $ n \ge 0 $, the integer $ \beta_n^R(M) := \rank_k\left( \Ext_R^n(M,k) \right) $ is called the $n$th {\it Betti number} of $M$. It is equal to the rank of the $ n $th component in a minimal free resolution of $M$. In other words, $ \beta_n^R(M) = \mu\left( \Omega_n^R(M) \right) $. Dually, $ \mu_n^R(M) := \rank_k\left( \Ext_R^n(k,M) \right) $ is called the $n$th {\it Bass number} of $M$. The formal sums $ P_M^R(t) := \sum_{n \ge 0} \beta_n^R(M) t^n $ and $ B_M^R(t) := \sum_{n \ge 0} \mu_n^R(M) t^n $ are called the {\it Poincar\'{e} series} and {\it Bass series} of $ M $ respectively. The following notions were introduced by Avramov.
 
 \begin{definition} (1) \cite[(3.1)]{Avr89}\label{def:complexity}
 	The {\it complexity} of $ M $, denoted $ \cx_R(M) $, is the smallest non-negative integer $ b $ such that $ \beta_n^R(M) \le \alpha n^{b-1} $ for all $ n \gg 0 $, and for some real number $ \alpha > 0 $. If no such $ b $ exists, then we set $ \cx_R(M) = \infty $.
 	
 	(2) \cite[p.~319]{Avr96} Replacing $ \beta_n^R(M) $ by $ \mu_n^R(M) $ in (1), one obtains the notion of {\it injective complexity} $ \injcx_R(M) $. This is called {\it plexity} $ \px_R(M) $ in \cite[(5.1)]{Avr89}.
 	
 	(3) \cite[p.~320]{Avr96} The {\it curvature} of $ M $, denoted $ \curv_R(M) $, is the reciprocal value of the radius of convergence of $ P_M^R(t) $, i.e.,
 	\[
 		\curv_R(M) := \limsup_{n \to \infty} \sqrt[n]{\beta_n^R(M)}.
 	\]
 \end{definition}

 We use the following elementary properties of complexity and curvature.

 \begin{lemma}\label{lem:cx-properties}
 	Let $ M $ and $ N $ be $ R $-modules. The following statements hold true.
 	\begin{enumerate}[\rm (i)]
 		\item
 		\cite[Rmk.~4.2.3]{Avr98}\\
 		{\rm (a)} $ \projdim_R(M) < \infty \; \Longleftrightarrow \; \cx_R(M) = 0 \; \Longleftrightarrow \; \curv_R(M) = 0 $.\\
 		{\rm (b)} $ \cx_R(M) < \infty \; \Longrightarrow \; \curv_R(M) \le 1 $; the converse is also true for $ M = k $ {\rm (}see \cite[8.1.2 and 8.2.2]{Avr98}{\rm )}.\\ 
 		{\rm (c)} $ \curv_R(M) $ is always finite.
 		\item
 		\cite[Prop.~4.2.4(2) and (3)]{Avr98}\\
 		{\rm (a)} $ \cx_R(M) = \cx_R\left( \Omega_n^R(M) \right) $ and $ \curv_R(M) = \curv_R\left( \Omega_n^R(M) \right) $ for all $ n \ge 1 $.\\
 		{\rm (b)} $ \cx_R(M\oplus N) = \max\{ \cx_R(M), \cx_R(N) \} $.\\
 		{\rm (c)} $ \curv_R(M\oplus N) = \max\{ \curv_R(M), \curv_R(N) \} $.	
 		\item
 		Let $ x \in R $ be regular on both $ R $ and $ M $. Then 
 		\[
 			\cx_{R/(x)}(M/xM) = \cx_R(M), \quad  \curv_{R/(x)}(M/xM) = \curv_R(M) 
 		\]
 		and \;$ \injcx_{R/(x)}(M/xM) = \injcx_R(M) $.
 	\end{enumerate}
 \end{lemma}
 
 \begin{proof}
 	(iii) Since $ x $ is regular on both $ R $ and $ M $, for every $ n \ge 0 $, we have that
 	\[
	 	\beta_n^{R/(x)}(M/xM) = \beta_n^R(M) \quad \mbox{and} \quad \mu^{R/(x)}_n(M/xM) = \mu^R_{n+1}(M);
 	\]
 	see, e.g., \cite[p.~140, Lem.~2]{Mat86}. Hence the desired equalities follow.
 \end{proof}  
 
 The following lemma shows that complexity and curvature of the residue field remain same after going modulo a regular element.
 
 \begin{lemma}\label{lem:cx-mod-x}
 	Let $ x \in \fm \smallsetminus \fm^2 $ be an $ R $-regular element. Then
 	\[
 		\cx_{R/(x)}(k) = \cx_R(k) \quad \mbox{and} \quad \curv_{R/(x)}(k) = \curv_R(k).
 	\]
 \end{lemma}

\begin{proof}
	Set $\overline{(-)} := (-) \otimes_R R/(x)$. In view of Lemma~\ref{lem:cx-properties}(ii) and (iii), we obtain that
	\begin{align*}
	\cx_R(k) & = \cx_R\left( \Omega_1^R(k) \right) = \cx_{\oR}\left( \overline{\Omega_1^R(k)} \right) \\
	& = \cx_{\oR}\left( \Omega_1^{\oR}(k) \oplus \Omega_0^{\oR}(k) \right) \quad \mbox{[by \cite[Cor.~5.3]{Tak06}]} \\
	& = \max\left\{ \cx_{\oR}\left( \Omega_1^{\oR}(k) \right), \cx_{\oR}\left( \Omega_0^{\oR}(k) \right) \right\} \\
	& = \cx_{\oR}(k).
	\end{align*}
	Similarly, one obtains that $ \curv_R(k) = \curv_{R/(x)}(k) $.
\end{proof}

 Let us recall a few well-known properties of complexities of the residue field.
 
 \begin{proposition}\label{prop:facts-projcx-k}{~}
 	\begin{enumerate}[\rm (i)]
 		\item
 		{\rm \cite[Prop.~2]{Avr96}} The residue field has maximal complexities and curvature:
 		\begin{align*}
 		\cx_R(k) & = \sup\{ \cx_R(M) : M \mbox{ is an $ R $-module} \} \\
 		& = \sup\{ \injcx_R(M) : M \mbox{ is an $ R $-module} \} = \injcx_R(k),\\
 		\mbox{and } \; \curv_R(k) & = \sup\{ \curv_R(M) : M \mbox{ is an $ R $-module} \}.
 		\end{align*}
 		\item
 		{\rm \cite[(2.3)]{Gul80}} If $ \cx_R(k) $ is finite, then $ R $ is a complete intersection ring.
 		\item 
 		If $ R $ is a complete intersection ring of codimension $c$, then it follows from \cite[Thm.~6]{Tat57} that $ \cx_R(k) = c $.
 		\item 
 		The statements {\rm (ii)} and {\rm (iii)} along with {\rm Lemma~\ref{lem:cx-properties}(i)(b)} yield the following:
 		\[
 			R \mbox{ is complete intersection } \; \Longleftrightarrow \; \cx_R(k) < \infty \; \Longleftrightarrow \; \curv_R(k) \le 1.
 		\]
 	\end{enumerate}
 \end{proposition}
 
 \begin{definition}\label{defn:max-complexity-curvature}
 	 In view of Proposition~\ref{prop:facts-projcx-k}(i), an $ R $-module $ M $ is said to have maximal complexity (resp., curvature) if $ \cx_R(M) = \cx_R(k) $ (resp., $ \curv_R(M) = \curv_R(k) $).
 \end{definition}
 
  We need the following elementary fact on vanishing of Exts or Tors.
  
 \begin{lemma}\cite[2.6]{Gho}\label{lem:Ext-Tor-mod-x}
 	Let $ M $ and $ N $ be $ R $-modules. Let $ x \in R $ be an $ R \oplus M \oplus N$-regular element. Set $\overline{(-)} := (-) \otimes_R R/(x)$. Suppose $ m $ and $ n $ are positive integers. If $ \Ext_R^i(M,N) = 0 $ {\rm(}resp. $ \Tor_i^R(M,N) = 0 ${\rm )} for all $ n \le i \le n + m $, then
 	\begin{align*}
 	&\Ext_{\oR}^i(\oM,\oN) = 0 \quad\mbox{for all } n \le i \le n + m - 1\\
 	(\mbox{resp.} \quad &\Tor_i^{\oR}(\oM,\oN) = 0 \quad\mbox{for all } n + 1 \le i \le n + m).
 	\end{align*}
 \end{lemma}
 
 The following proposition is a consequence of a result due to Huneke, \c{S}ega and Vraciu \cite[Prop.~2.9]{HSV04}.
 
 \begin{proposition}\label{prop:HSV04-2.9}
 	Let $ R $ be a non-Gorenstein local ring such that $ \fm^3 = 0 $. Let $ M $ be an $ R $-module with the property that $ \Tor_i^R(M,E) = 0 $ for all $ i \gg 0 $, where $ E := E_R(k) $. Then $ \beta_n^R(M) = c $ for every $ n \ge 1 $, where $ c $ is a constant.
 \end{proposition}
 
 \begin{proof}
 	Let us denote $ \Omega_1^R(M) $ by $ \Omega(M) $. If $ M $ is free, then $ \beta_n^R(M) = 0 $ for every $ n \ge 1 $. So we may assume that $ M $ is not free. Hence, by the Auslander-Buchsbaum Formula, $ \Omega(M) $ is also not free. Note that
 	\[
	 	\Tor_i^R(\Omega(M),E) \cong \Tor_{i+1}^R(M,E) = 0 \quad \mbox{for all } i \gg 0.
	\]
	Since $ \fm^3 = 0 $, we have that $ \fm^2 \subseteq \Soc(R) \subseteq \ann_R( \Omega(M) ) $ (by Lemma~\ref{lemma: socle, syzygy}), and hence $ \fm^2 \Omega(M) = 0 $. Therefore, by virtue of \cite[Prop.~2.9]{HSV04}, one can deduce that $ \beta_n^R(\Omega(M)) = c $ for every $ n \ge 0 $, where $ c $ is a nonzero constant. Hence $ \beta_n^R(M) = \beta_{n-1}^R(\Omega(M)) = c $ for every $ n \ge 1 $.
 \end{proof}
 
\section{Gorensteinness of short local rings}\label{sec:main}
 
 The following theorem shows that modules with extremal complexity or curvature can be used to detect whether a short local ring is Gorenstein.
 
 \begin{theorem}\label{thm:Gor-rings-max-cx}
 	Let $ (R,\mathfrak{m},k) $ be a short local ring. Let $ M $ be an $ R $-module such that either $ \curv_R(M) = \curv_R(k) $ or $ \cx_R(M) = \cx_R(k) $. {\rm (}Particularly, $ \cx_R(M) $ can be infinite{\rm )}. Then the following statements are equivalent:
 	\begin{enumerate}
 		\item[\rm (i)] $ R $ is Gorenstein.
 		\item[\rm (ii)] $ \Ext_R^i(M,R) = 0 $ for all $ i \gg 0 $.
 	\end{enumerate}
 	If $ R $ has a canonical module $ \omega $, then we may add the following:
 	\begin{enumerate}
 		\item[\rm (iii)] $ \Tor_i^R(M, \omega) = 0 $ for all $ i \gg 0 $.
 	\end{enumerate}
 \end{theorem}
 
 \begin{proof}
 	Since $ \underline{x} = x_1,\ldots,x_d $ is an $ R $-regular sequence such that $ \mathfrak{m}^3 \subseteq (\underline{x}) $, it can be easily verified that $ R $ is Cohen-Macaulay of dimension $ d $. If $ R $ is Gorenstein, then $ \Ext_R^i(M,R) = 0 $ for all $ i \ge d + 1 $ (because in this case $ \injdim_R(R) = d $). Moreover, if $ R $ is Gorenstein, then $ \omega \cong R $, and hence $ \Tor_i^R(M, \omega) = 0 $ for all $ i \ge 1 $. So we only need to prove the implications \{(ii) $ \Rightarrow $ (i)\} and \{(iii) $ \Rightarrow $ (i)\}.
 	
 	\{(ii) or (iii)\} $ \Rightarrow $ (i): We consider the case when $ \cx_R(M) = \cx_R(k) $. Another case, i.e., $ \curv_R(M) = \curv_R(k) $ can be treated similarly. We prove these implications by using induction on $ d $. Assume that $ d = 0 $. In this case, $ \fm^3 = 0 $, and the injective hull $ E $ (of $ k $ over $ R $) is a canonical module of $ R $, i.e., $ \omega \cong E $. Note that $ \Ext_R^i(M,R)^{\vee} \cong \Tor_i^R(M,E) $ for every $ i \ge 0 $. Hence, by Matlis Duality, $ \Ext_R^i(M,R) = 0 $ \iff $ \Tor_i^R(M,E) = 0 $. So, in this case, (ii) and (iii) are equivalent. We assume that $ \Tor_i^R(M,E) = 0 $ for all $ i \gg 0 $. If possible, assume that $ R $ is not Gorenstein. Then, by virtue of Proposition~\ref{prop:HSV04-2.9}, we obtain that $ \beta_n^R(M) = c $ for every $ n \ge 1 $, where $ c $ is a constant. Therefore $ \cx_R(k) = \cx_R(M) \le 1 $. Hence, in view of Proposition~\ref{prop:facts-projcx-k}(iv), $ R $ is a complete intersection ring, which contradicts the assumption that $ R $ is not Gorenstein. Therefore $ R $ is Gorenstein.
 	
 	We now give the inductive step. Assume that $ d \ge 1 $. It is given that $ \mathfrak{m}^3 \subseteq (\underline{x}) $, where $ \underline{x} = x_1,\ldots,x_d \in \fm \smallsetminus \fm^2 $ is $ R $-regular.
 	We set $ \overline{(-)} := (-) \otimes_R R/(x_1) $. It can be observed that $ (\fm/(x_1))^3 \subseteq (\overline{x_2},\ldots,\overline{x_d}) $, where $ \overline{x_2},\ldots,\overline{x_d} $ (the images of $ x_2,\ldots,x_d $ in $ \oR $ respectively) is an $ \oR $-regular sequence.
 	Suppose $ \Ext_R^i(M,R) = 0 $ (resp., $ \Tor_i^R(M, \omega) = 0 $) for all $ i \gg 0 $. Then
 	\begin{align}
	 	\Ext_R^i(\Omega(M), R) & \cong \Ext_R^{i+1}(M,R) = 0 \quad \mbox{for all } i \gg 0 \label{eqn:thm:Gor-rings-max-cx-2}\\
	 	(\mbox{resp., } \Tor_i^R(\Omega(M), \omega) & \cong \Tor_{i+1}^R(M, \omega) = 0 \quad \mbox{for all } i \gg 0).\nonumber
 	\end{align}
 	Since $ x_1 $ is $ R $-regular and $ \omega $ is a maximal Cohen-Macaulay $ R $-module, $ x_1 $ is regular on $ \omega $. Note that $ x_1 $ is also $ \Omega(M) $-regular. Let us denote $ \Omega(M) $ by $ N $. Hence, in view of \eqref{eqn:thm:Gor-rings-max-cx-2}, by Lemma~\ref{lem:Ext-Tor-mod-x}, we can deduce that $\Ext_{\oR}^i(\oN, \oR) = 0 $ (resp., $ \Tor_i^{\oR}(\oN, \oomega) = 0 $) for all $ i \gg 0 $. (Note that $ \oomega $ is a canonical module of $ \oR $.) On the other hand, by using Lemmas~\ref{lem:cx-properties} and \ref{lem:cx-mod-x}, we obtain that
 	\[
	 	\cx_{\oR}(\oN) = \cx_R(N) = \cx_R(M) = \cx_R(k) = \cx_{\oR}(k).
 	\]
 	Thus there is an $ \oR $-module $ \oN $ such that $ \cx_{\oR}(\oN) = \cx_{\oR}(k) $ and $\Ext_{\oR}^i(\oN, \oR) = 0 $ (resp., $ \Tor_i^{\oR}(\oN, \oomega) = 0 $) for all $ i \gg 0 $. Since $ \dim(\oR) = d - 1 $, in either case, by induction hypothesis, we get that $ \oR $ is Gorenstein, and hence $ R $ is Gorenstein.
 \end{proof}

\begin{remark}\label{rmk:thm-max-curv-cx}
	The two conditions `$ \curv_R(M) = \curv_R(k) $' and `$ \cx_R(M) = \cx_R(k) $' in Theorem~\ref{thm:Gor-rings-max-cx} are independent of one another. In Example~\ref{exam1:thm:char-of-reg}, assuming $ p > q $, we have $ \curv_R(N) < \curv_R(k) $, but $ \cx_R(N) = \infty = \cx_R(k) $. On the other hand, in Example~\ref{exam2:thm:char-of-reg}, $ \curv_R(N) = \curv_R(k) $. Since $ \beta_n^R(N) = 1 $ for all $ n \ge 0 $, we have $ \cx_R(N) = 1 < 2 = \cx_R(k) $ (by Proposition~\ref{prop:facts-projcx-k}(iii)).
\end{remark}

 As an immediate corollary of Theorem~\ref{thm:Gor-rings-max-cx}, we obtain one of our main results.
 
 \begin{corollary}\label{cor:Gor-rings-hom-image-of-syz}
 	Let $ (R,\mathfrak{m},k) $ be a short local ring. Let $ \omega $ be a canonical module of $ R $. Let $ M $ be a nonzero homomorphic image of a finite direct sum of syzygy modules of $ k $. Then the following statements are equivalent:
 	\begin{enumerate}[\rm (i)]
 		\item $ R $ is Gorenstein.
 		\item $ \Ext_R^i(M,R) = 0 $ for all $ i \gg 0 $.
 		\item $ \Tor_i^R(\omega, M) = 0 $ for all $ i \gg 0 $.
 	\end{enumerate}
 \end{corollary}
 
 \begin{proof}
 	The corollary follows from Theorems~\ref{thm:Avr-max-cx} and \ref{thm:Gor-rings-max-cx}.
 \end{proof}

 As some other consequences  of Theorem~\ref{thm:Gor-rings-max-cx}, we obtain the following results.

\begin{para}\label{para2:thm-max-curv-cx}
	With $ R $ and $ \omega $ as in Corollary~\ref{cor:Gor-rings-hom-image-of-syz}, we may recover that if $ \Ext_R^i(\omega, R) = 0 $ (resp., $ \Tor_i^R(\omega, \omega) = 0 $) for all $ i \gg 0 $, then $ R $ is Gorenstein.
\end{para}

\begin{proof}
	Let $ \Ext_R^i(\omega, R) = 0 $ (resp., $ \Tor_i^R(\omega, \omega) = 0 $) for all $ i \gg 0 $. If possible, suppose that $ R $ is not Gorenstein. Then $ \cx_R(\omega) = \infty $ by \cite[Prop.~1.1.(4)]{JL07}. Hence Theorem~\ref{thm:Gor-rings-max-cx} yields that $ R $ is Gorenstein, a contradiction.
\end{proof}

\begin{remark}
	In \cite[Thm.~2.11]{HSV04}, it is shown that if $ \fm^3 = 0 $ and $ \Ext_R^i(\omega, R) = 0 $ for any three consecutive $ i \ge 0 $, then $ R $ is Gorenstein.
\end{remark}
 \begin{para}\label{para23:thm-max-curv-cx}
 	Let $ Q = k[X_1,\ldots,X_e] $ be a positively graded polynomial ring over a field $ k $ of characteristic $ 0 $, and $ I $ be a homogeneous ideal containing $ (X_1,\ldots,X_e)^3 $. Set $ R := Q/I $. Let $ C $ denote either $ I/I^2 $ or the module $ \Omega_{R| k} $ of K\"{a}hler differentials over $ k $. If $ \Ext_R^i(C, R) = 0 $ (resp., $ \Tor_i^R(C, \omega) = 0 $) for all $ i \gg 0 $, then $ R $ is Gorenstein.
 \end{para}

\begin{proof}
	If $ R $ is not Gorenstein, then by \cite[Thm.~8.3.3 and 8.3.4]{Avr98}, $ \curv_R(C) = \curv_R(k) $. Hence the proof is similar as that of \ref{para2:thm-max-curv-cx}.
\end{proof}

 The following result is a counterpart of Theorem~\ref{thm:Gor-rings-max-cx} for injective complexity.
 
 \begin{theorem}\label{thm:Gor-rings-max-injcx}
 	With $ R $ and $ \omega $ as in {\rm Corollary~\ref{cor:Gor-rings-hom-image-of-syz}}, let $ M $ be an $ R $-module such that $ \injcx_R(M) = \injcx_R(k) $ and $ \Ext_R^i(\omega, M) = 0 $ for all $ i \gg 0 $. Then $ R $ is Gorenstein.
 \end{theorem}
 
 \begin{proof}
 	Let $ \Ext_R^i(\omega, M) = 0 $ for all $ i \gg 0 $. Note that $ \underline{x} = x_1,\ldots,x_d $. We use induction on $ d $. In the base case, i.e., if $ d = 0 $, then $ \omega \cong E $ $ ( := E_R(k)) $. Hence $ \Tor_i^R(E,M^{\vee}) \cong \Ext_R^i(E,M)^{\vee} = 0 $ for all $ i \gg 0 $. Note that
 	\begin{align}
	 	\beta_n^R(M^{\vee}) & = \rank_k\left( \Tor_n^R(k, M^{\vee}) \right) = \rank_k\left( \Ext_R^n(k,M)^{\vee} \right) \label{eqn:thm:Gor-rings-max-injcx-0}\\
	 	& = \rank_k( \Ext_R^n(k,M) ) \quad \mbox{[as $ k^{\vee} \cong k $]} \nonumber \\
	 	& = \mu_n^R(M) \quad \mbox{for every $ n \ge 0 $}. \nonumber 
 	\end{align}
 	Thus $ \cx_R(M^{\vee}) = \injcx_R(M) = \injcx_R(k) = \cx_R(k) $. So, in view of the implication \{(iii) $ \Rightarrow $ (i)\} in Theorem~\ref{thm:Gor-rings-max-cx}, we obtain that $ R $ is Gorenstein.
 	
 	We now give the inductive step. Assume that $ d \ge 1 $. 
 	Set $ \overline{(-)} := (-) \otimes_R R/(x_1) $. By virtue of \cite[Thm.~A]{AB89}, we have a maximal Cohen-Macaulay (MCM) approximation of $ M $, i.e., a short exact sequence $ 0 \to Y \to N \to M \to 0 $ of $ R $-modules, where $ N $ is MCM and $ Y $ has finite injective dimension. Hence
 	\begin{equation}\label{eqn:thm:Gor-rings-max-injcx-1}
	 	 \Ext_R^i(\omega,N) \cong \Ext_R^i(\omega,M) = 0 \quad \mbox{and} \quad \Ext_R^i(k,N) \cong \Ext_R^i(k,M)
 	\end{equation}
 	for all $ i \gg 0 $. In particular, $ \mu_n^R(N) = \mu_n^R(M) $ for all $ n \gg 0 $, which yields that $ \injcx_R(N) = \injcx_R(M) = \cx_R(k) $. Since $ \omega $ and $ N $ are MCM $ R $-modules and $ x_1 $ is $ R $-regular, it follows that $ x_1 $ is regular on both $ \omega $ and $ N $. So, in view of \eqref{eqn:thm:Gor-rings-max-injcx-1}, by Lemma~\ref{lem:Ext-Tor-mod-x}, one obtains that $\Ext_{\oR}^i\big( \overline{\omega}, \oN \big) = 0 $ for all $ i \gg 0 $. It follows from Lemma~\ref{lem:cx-properties}(iii) and Lemma~\ref{lem:cx-mod-x} that $	\injcx_{\oR}(\oN) = \injcx_R(N) = \cx_R(k) = \cx_{\oR}(k) $. Thus there is an $ \oR $-module $ \oN $ such that $	\injcx_{\oR}(\oN) = \cx_{\oR}(k) $ and $\Ext_{\oR}^i\big( \overline{\omega}, \oN \big) = 0 $ for all $ i \gg 0 $. Therefore, since $ \dim(\oR) = d - 1 $, by induction hypothesis, we get that $ \oR $ is Gorenstein, and hence $ R $ is Gorenstein.
 \end{proof}
 
 As a few consequences of Theorem~\ref{thm:Gor-rings-max-injcx}, we obtain the following:
 
 \begin{corollary}\label{cor:Gor-Artin-rings-mM-and-M-mod-colon}
 	With $ R $ and $ \omega $ as in {\rm Corollary~\ref{cor:Gor-rings-hom-image-of-syz}}, let $ M $ be an $ R $-module.
 	\begin{enumerate}[\rm (i)]
 		\item If $ \fm M \neq 0 $ and $ \Ext_R^i(\omega,\fm M) = 0 $ for all $ i \gg 0 $, then $ R $ is Gorenstein.
 		\item Let $ n \ge 1 $, such that $ (0 :_M \fm^n) \neq (0 :_M \fm^{n+1}) $. If $ \Ext_R^i(\omega, M/(0 :_M \fm^n)) = 0 $ for all $ i \gg 0 $, then $ R $ is Gorenstein.
 	\end{enumerate}
 \end{corollary}
  
 \begin{proof}
 	With the above hypotheses, by \cite[Thm.~4 and Prop.~7]{Avr96},
 	\[
 		\injcx_R(\fm M) = \injcx_R(k) \quad \mbox{ and} \quad \injcx_R( M/(0 :_M \fm^n) ) = \injcx_R(k).
 	\]
 	Hence the proof follows from Theorem~\ref{thm:Gor-rings-max-injcx}.
 \end{proof}
 
 \begin{remark}\label{rmk:cor:Gor-Artin-rings-mM-and-M-mod-colon}
 	Let $ M $ be a module over an arbitrary local ring $ R $ such that $ \fm M \neq 0 $. It is true in general that if $ \Ext_R^i(\omega,\fm M) = 0 $ for all $ i \gg 0 $, then $ \projdim_R(\omega) $ is finite, and hence $ R $ is Gorenstein; see \cite[Thm.~1.5(a)]{TTY07} and the proof of \cite[Thm.~1.1]{LV68}. But Corollary~\ref{cor:Gor-Artin-rings-mM-and-M-mod-colon}(ii) was not known. However, we would like to state Corollary~\ref{cor:Gor-Artin-rings-mM-and-M-mod-colon} as applications of Theorem~\ref{thm:Gor-rings-max-injcx}.
 \end{remark}

\begin{corollary}\label{cor:Gupta}
	With $ R $ and $ \omega $ as in {\rm Corollary~\ref{cor:Gor-rings-hom-image-of-syz}}, let $ M $ be a nonzero homomorphic image of a finite direct sum of copies of $ \Omega_1^R(k) $ $( = \fm )$ and $ \Omega_0^R(k) $ $( = k)$. Then $ R $ is Gorenstein \iff $ \Ext_R^i(\omega, M) = 0 $ for all $ i \gg 0 $.
\end{corollary}

\begin{proof}
	We only need to prove the `if' part. Let $ \Ext_R^i(\omega, M) = 0 $ for all $ i \gg 0 $. In view of Lemma~\ref{lem:Gupta}, $ M \cong \fm L \oplus k^{\oplus u} $ for some $ R $-module $ L $ and integer $ u \ge 0 $. Hence $ \Ext_R^i(\omega, \fm L) = 0 $ for all $ i \gg 0 $. Therefore, if $ \fm L \neq 0 $, by Corollary~\ref{cor:Gor-Artin-rings-mM-and-M-mod-colon}(i), $ R $ is Gorenstein. So we may assume that $ \fm L = 0 $, hence $ u \ge 1 $. Thus $ \Ext_R^i(\omega, k) = 0 $ for all $ i \gg 0 $, which yields that $ \pd_R(\omega) < \infty $, equivalently, $ R $ is Gorenstein.
\end{proof}

\begin{remark}
	Corollary~\ref{cor:Gupta} provides a partial answer to Question~\ref{ques2}.
\end{remark}

  The author is grateful to Anjan Gupta for making the following:
  
\begin{lemma}\label{lem:Gupta}
	Over an arbitrary local ring $ (R,\fm,k) $, a homomorphic image of a finite direct sum of copies of $ \Omega_1^R(k) $ $( = \fm )$ and $ \Omega_0^R(k) $ $( = k )$ can be written as $ \fm L \oplus k^{\oplus u} $ for some $ R $-module $ L $ and integer $ u \ge 0 $, where $ k^{\oplus u} $ denotes the direct sum of $ u $ many copies of $ k $.
\end{lemma}

\begin{proof}
	Let $ M $ be an $ R $-module, and $ N $ be a submodule of $ M \oplus k $. We claim that
	\begin{equation}\label{prop:Gupta-eqn0}
		\mbox{$ (M \oplus k)/N $\; is isomorphic to either\; $ (M/N) \oplus k $\; or\; $ M/M' $}
	\end{equation}
	for some submodule $ M' $ of $ M $. To prove the claim, assume that $ N $ is generated by $ (x_1,\overline{a_1}), \ldots, (x_n,\overline{a_n}) $ for some $ x_i \in M $ and $ a_i \in R $. If $ \overline{a_i} = 0 $ for all $ 1 \le i \le n $, then $ N $ is a submodule of $ M $, and hence $ (M \oplus k)/N \cong (M/N) \oplus k $. In another case, without loss of generality, we may assume that $ a_1 = 1 $. Note that
	\begin{equation}\label{prop:Gupta-eqn1}
		\dfrac{M \oplus k}{\langle (x_1,\overline{1}) \rangle} \cong \dfrac{(M \oplus R)/(0 \oplus \fm)}{\langle (x_1,\overline{1}) \rangle} \cong \dfrac{M \oplus R}{\langle (x_1,1) \rangle + (0 \oplus \fm)}.
	\end{equation}
	The map $ \varphi : M \oplus R \to M $ given by $ \varphi((y,b)) = y - b x_1 $ induces an isomorphism
	\begin{equation}\label{prop:Gupta-eqn2}
	(M \oplus R)/\langle (x_1,1) \rangle \cong M.
	\end{equation}
	Thus we have
	\begin{align*}
	\dfrac{M \oplus k}{N} & \cong \dfrac{(M \oplus k)/\langle (x_1,\overline{1}) \rangle}{N'} \quad \mbox{for some } N' \\
	& \cong \dfrac{(M \oplus R)/\big( \langle (x_1,1) \rangle + (0 \oplus \fm) \big)}{N'}  \quad \mbox{[by \eqref{prop:Gupta-eqn1}]} \\
	& \cong \dfrac{(M \oplus R)/\langle (x_1,1) \rangle}{M'}  \quad \mbox{for some } M' \\
	& \cong M/M' \quad \mbox{[by \eqref{prop:Gupta-eqn2}]}.
	\end{align*}
	Using \eqref{prop:Gupta-eqn0} repeatedly, one can deduce that
	\[
		( \fm^{\oplus s} \oplus k^{\oplus t} ) /U \cong \big(( \fm^{\oplus s} )/U'\big) \oplus k^{\oplus u} \cong \fm (R^s/U') \oplus k^{\oplus u}
	\]
	for some $ 0 \le u \le t $ and submodule $ U' $ of $ \fm^{\oplus s} $. Now set $ L := R^s/U' $ to get the desired result.
\end{proof}

With hypotheses as in Corollary~\ref{cor:Gor-rings-hom-image-of-syz}, we now investigate whether $ R $ is Gorenstein when $ \Ext_R^i(\omega,M) = 0 $ for all $ i \gg 0 $. Let us consider the case when $ R $ is Artinian.

 \begin{proposition}\label{prop:length-neq-2e-Vanishing-of-Ext}
 	Suppose $ \fm^3 = 0 $, and $ 2 \cdot \mu(\fm) \neq \lambda(R) $. Set $ E := E_R(k) $. Let $ M $ be an $ R $-module such that $ \Ext_R^i(E,M) = 0 $ for three consecutive values of $ i \ge 3 $. Then either $ E \cong R $ or $ M $ is injective.
 \end{proposition}
 
 \begin{proof}
 	If possible, assume that $ E \ncong R $ and $ M $ is not injective. Note that $ E \ncong R $ is equivalent to that $ R $ is not Gorenstein. Moreover, since $ E $ is indecomposable as an $ R $-module, $ E $ is not free. As $ M $ is not injective, by Matlis Duality, we can deduce that $ M^{\vee} $ is not free. Hence it follows that $ \Omega(M^{\vee}) $ is not free. Since $ \fm^3 = 0 $, we have $ \fm^2 \Omega(M^{\vee}) = 0 $ (by Lemma~\ref{lemma: socle, syzygy}). Note that
 	\[
	 	\Tor_j^R(E,\Omega(M^{\vee})) \cong \Tor_{j+1}^R(E,M^{\vee}) \cong \Ext_R^{j+1}(E,M)^{\vee} = 0
 	\]
 	for three consecutive values of $ j \ge 2 $. Therefore, in view of \cite[Rmk.~2.4]{HSV04}, we have $ \Soc(R) = \fm^2 $. Moreover, by virtue of \cite[Prop.~2.9]{HSV04}, we get that $ e = t + 1 $, where $ e := \mu(\fm) $ and $ t := \rank_k(\Soc(R)) $. Hence $ \lambda(R) = \lambda(R/\fm) + \lambda(\fm/\fm^2) + \lambda(\fm^2) = 1 + e + t = 2e $, which is a contradiction. So either $ E \cong R $ or $ M $ is injective.
 \end{proof}
 
 Without the condition `$ 2 \cdot \mu(\fm) \neq \lambda(R) $', Proposition~\ref{prop:length-neq-2e-Vanishing-of-Ext} is not necessarily true.
 
 \begin{example}\label{exam:cannot-omit-length-neq-2e}
 	Let $ R = k[x,y,z]/(x^2,xy,y^2,z^2) $, where $ k $ is a field. We collect this ring from \cite[Example~2.8]{JL07}. Note that $ \fm^3 = 0 $, where $ \fm := (x, y, z) $ is the maximal ideal of $ R $. Moreover, $ \mu(\fm) = 3 $ and $ \lambda(R) = 6 $. Therefore $ 2 \cdot \mu(\fm) = \lambda(R) $. Since $ \Soc(R) = (xz,yz) $, it follows that $ R $ is not Gorenstein, hence $ E \ncong R $. We set $ N := (z) $, and $ M := N^{\vee} $. Since $ N $ is annihilated by $ z $, it is not free. Hence we can deduce that $ M $ is not injective. By virtue of Matlis Duality, we obtain that
 	\begin{align}
	 	\Ext_R^i(E,M) & \cong \Ext_R^i(E,M)^{\vee \vee} \cong \Tor_i^R(E,M^{\vee})^{\vee} \cong \Tor_i^R(E,N)^{\vee} \label{Ext-Tor-iso}\\
	 	& \cong \Tor_i^R(N, E)^{\vee} \cong \Ext_R^i(N, R) = 0 \quad \mbox{for every } i \ge 1.\nonumber
 	\end{align}
 	To get the last equality, one may compute $ \Ext_R^i(N, R) $ by considering the minimal free resolution of $ N $:\quad $ \cdots \stackrel{z}{\longrightarrow} R \stackrel{z}{\longrightarrow} R \stackrel{z}{\longrightarrow} R \to 0 $.
 \end{example}

\begin{remark}
	In Example~\ref{exam:cannot-omit-length-neq-2e}, we should note that $ \Ext_R^i(E,M) = 0 $ for every $ i \ge 1 $, but neither of $ \projdim_R(E) $ and $ \injdim_R(M) $ is finite. Moreover, $ \Tor_i^R(N, E) = 0 $ for every $ i \ge 1 $, but neither of $ \projdim_R(N) $ and $ \projdim_R(E) $ is finite.
\end{remark}

 As a consequence of Proposition~\ref{prop:length-neq-2e-Vanishing-of-Ext}, we obtain the following:
 
 \begin{corollary}\label{cor:length-neq-2e-Gor-rings}
 	Let $ \fm^3 = 0 $, and $ 2 \cdot \mu(\fm) \neq \lambda(R) $. Set $ E := E_R(k) $. Let $ M $ be a nonzero homomorphic image of a finite direct sum of syzygy modules of $ k $. If $ \Ext_R^i(E,M) = 0 $ for three consecutive values of $ i \ge 3 $, then $ R $ is Gorenstein.
 \end{corollary}
 
 \begin{proof}
 	If $ \Ext_R^i(E,M) = 0 $ for three consecutive values of $ i \ge 3 $, then by Proposition~\ref{prop:length-neq-2e-Vanishing-of-Ext}, either $ E \cong R $ or $ M $ is injective. If $ E \cong R $, then $ R $ is Gorenstein. If $ M $ is injective, then in view of \cite[Cor.~3.4]{GGP}, we obtain that $ R $ is regular, and hence $ R $ is Gorenstein. Thus, in both cases, $ R $ is Gorenstein.
 \end{proof}
 
 The following proposition provides us another class of modules $ M $ for which the vanishing of $ \Ext_R^{\gg 0}(E,M) $ ensures that $ R $ is Gorenstein.
 
 \begin{proposition}\label{prop:mu-le-length-mM-vanishing-of-Ext}
 	Let $ \fm^3 = 0 $. Let $ M $ be a nonzero $ R $-module such that $ \fm^2 M = 0 $ and $ \mu(M) \le \lambda(\fm M) $. If $ \Ext_R^i(E,M) = 0 $ for all $ i \gg 0 $, then $ R $ is Gorenstein.
 \end{proposition}
 
 \begin{proof}
 	Let $ \Ext_R^i(E,M) = 0 $ for all $ i \gg 0 $. Hence $ \Tor_i^R(E,M^{\vee}) = 0 $ for all $ i \gg 0 $. If possible, assume that $ R $ is not Gorenstein. Then, by Proposition~\ref{prop:HSV04-2.9}, we have $ \beta_n^R(M^{\vee}) = c $ for every $ n \ge 1 $, where $ c $ is a constant. Thus, in view of \eqref{eqn:thm:Gor-rings-max-injcx-0}, we obtain that $ \mu_n^R(M) = c $ for every $ n \ge 1 $. Set $ a := \mu(M) $ and $ b := \lambda(\fm M) $. Note that $ \fm M $ is annihilated by $ \fm $. Moreover, $ \rank_k(\fm M) = b $ and $ \rank_k(M/\fm M) = a $. So there is a short exact sequence $ 0 \to k^b \to M \to k^a \to 0 $, which yields an exact sequence $ \Ext_R^{n - 1}(k,k^a) \to \Ext_R^n(k,k^b) \to \Ext_R^n(k,M) $ 	for every $ n \ge 1 $. From this exact sequence, for every $ n \ge 1 $, it follows that
 	\begin{align*}
	 	b \cdot \beta_n^R(k) & \le a \cdot \beta_{n-1}^R(k) + \mu_n^R(M) = a \cdot \beta_{n-1}^R(k) + c \\
	 	& \le b \cdot \beta_{n-1}^R(k) + c \quad \mbox{[because $ a \le b $]} \\
	 	& \le a \cdot \beta_{n-2}^R(k) + 2 c \quad \mbox{[using the 1st inequality for $ n - 1 $]} \\
	 	& \;\; \vdots \\
	 	& \le a \cdot \beta_0^R(k) + n c = c n + a.
 	\end{align*}
 	Therefore, since $ b \ge a > 0 $, we get that $ \beta_n^R(k) \le (c/b) \cdot n + (a/b) $ for every $ n \ge 1 $. This implies that $ \cx_R(k) $ is finite. Hence, by Proposition~\ref{prop:facts-projcx-k}(ii), $ R $ is a complete intersection ring, which contradicts the assumption that $ R $ is not Gorenstein. So $ R $ must be Gorenstein.
 \end{proof}
 
 \begin{remark}\label{rmk:prop:mu-le-length-mM-vanishing-of-Ext}
 	\begin{enumerate}[{\rm (i)}]
 		\item In Proposition~\ref{prop:mu-le-length-mM-vanishing-of-Ext}, $ M $ can be taken as $ R/I $, where $ I $ is an ideal of $ R $ such that $ \fm^2 \subseteq I \subsetneq \fm $.
 		\item Let $ M $ be a homomorphic image of a finite direct sum of syzygy modules of $ k $. We should note that if $ \fm^3 = 0 $, then we have $ \fm^2 M = 0 $, but $ M $ need not satisfy the condition `$ \mu(M) \le \lambda(\fm M) $'.
 	\end{enumerate}
 \end{remark}
 
 The condition `$ \mu(M) \le \lambda(\fm M) $' cannot be omitted from Proposition~\ref{prop:mu-le-length-mM-vanishing-of-Ext}.
 
 \begin{example}\cite[Example~2.8]{JL07}\label{exam:cannot-omit-mu-le-length-mM}
 	Let $ R $, $ M $ and $ N $ be as in Example~\ref{exam:cannot-omit-length-neq-2e}. Note that $ \fm^3 = 0 $, $ M \neq 0 $ and $ \fm^2 M = 0 $. Moreover, $ \Ext_R^i(E,M) = 0 $ for every $ i \ge 1 $, but $ R $ is not Gorenstein. One can verify that $ \mu(M) \nleqslant \lambda(\fm M) $. Indeed,
 	\begin{align}
	 	\mu(M) & = \type(N) \quad \mbox{[see, e.g., \cite[3.2.12(d)]{BH98}]} \label{mu-M-Soc-N}\\
	 	& = \rank_k((0 :_N \fm)) = \rank_k((xz,yz)) = 2 \nonumber
 	\end{align}
 	and $ \lambda(\fm M) = \lambda(M) - \mu(M) = \lambda(N) - 2 = 1 $.
 \end{example}
 
 As a consequence of Proposition~\ref{prop:mu-le-length-mM-vanishing-of-Ext}, we can recover a result of Jorgensen and Leuschke (which is a part of \cite[Thm.~2.5]{JL07}).
 
 \begin{corollary}\label{cor:lenth-mN-le-mu-N-vanishing-Ext}
 	Let $ \fm^3 = 0 $. Let $ N $ be a nonzero $ R $-module such that $ \fm^2 N = 0 $ and $ \lambda(\fm N) \le \mu(N) $. If $ \Ext_R^i(N,R) = 0 $ for all $ i \gg 0 $, then $ R $ is Gorenstein.
 \end{corollary}

\begin{proof}
	If $ k $ is a direct summand of $ N $, then $ \Ext_R^i(N,R) = 0 $ for some $ i \ge 1 $ implies that $ R $ is Gorenstein. So we may assume that $ k $ is not a direct summand of $ N $. Set $ M := N^{\vee} $. In view of \eqref{mu-M-Soc-N}, $ \mu(M) = \rank_k(\Soc(N)) = \lambda(\fm N) $ by \cite[2.3]{HSV04}. Hence $ \lambda(\fm M) = \lambda(M) - \mu(M) = \lambda(N) - \lambda(\fm N) = \mu(N) $. Thus $ \fm^2 M = 0 $ and $ \mu(M) \le \lambda(\fm M) $. In view of \eqref{Ext-Tor-iso}, $ \Ext_R^i(E,M) \cong \Ext_R^i(N, R) $ for every $ i \ge 1 $. So the result follows from Proposition~\ref{prop:mu-le-length-mM-vanishing-of-Ext}.
\end{proof}

  With the hypotheses as in Theorem~\ref{thm:Gor-rings-max-cx}, one may ask that if $ \Ext^i_R(\omega, M) = 0 $ for all $ i \gg 0 $, then is $ R $ Gorenstein? In this situation, $ R $ is not necessarily Gorenstein. The author is grateful to Ryo Takahashi for pointing out this fact with the following:
 
 \begin{example}\label{exam:Takahashi}
 	Consider a non-Gorenstein local ring $ (R,\fm,k) $ such that $ \fm^3 = 0 $. In this case, Betti number of the canonical module $ \omega $ grows exponentially, due to \cite[Prop.~1.1.(4)]{JL07}. Therefore $ \cx_R(\omega) = \infty $. In this situation, although $ \cx_R(\omega) = \cx_R(k) $ and $ \Ext^i_R(\omega, \omega) = 0 $ for all $ i \ge 1 $, but $ R $ is not Gorenstein. As an example, if $ R = k[x,y]/(x^2,xy,y^2) $, then $ \cx_R(\omega) = \cx_R(k) $. Moreover $ \curv_R(\omega) = \curv_R(k) $ (cf., \cite[4.2.2]{Avr98}) and $ \Ext^i_R(\omega, \omega) = 0 $ for all $ i \ge 1 $, but $ R $ is not Gorenstein.
 \end{example}

\section{Why is Gorensteinness of a short local ring important?}\label{sec:last}
 
 The class of Gorenstein short local rings is well studied. For instance, Poincar\'{e} series and Bass series of all (finitely generated) modules over such a ring are rational, sharing a common denominator; see \cite{Sjo79} and \cite{MS18}. The Koszulness of modules over Gorenstein short local rings are studied in \cite[Thm.~4.6 and Cor.~4.7]{AIS08}. Moreover, it is shown in \cite[Thm.~4.1]{AIS08} when such rings are Koszul. 
 So there are enough reasons to study Gorensteinness of a short local ring. In this theme, our contribution is the following theorem which provides a necessary and sufficient condition for a Gorenstein short local ring to be regular.
 
 \begin{theorem}\label{thm:char-of-reg}
 	Suppose $ R $ is a Gorenstein short local ring.	Let $ M $ and $ N $ be $ R $-modules having maximal complexity.	$ ( $Possibly, $ M = N ) $. Then $ R $ is regular if and only if $ \Tor_i^R(M,N) = 0 $ for all $ i \gg 0 $.
 \end{theorem}

\begin{proof}	
	We only need to prove the `if' part. Let $ \Tor_i^R(M,N) = 0 $ for all $ i \gg 0 $. Since $ \Tor_i^R(\Omega(M),\Omega(N)) \cong \Tor_{i+1}^R(M,\Omega(N)) \cong \Tor_{i+2}^R(M,N) $ for every $ i \ge 1 $, and $ \cx_R(M) = \cx_R(\Omega(M)) $, $ \cx_R(N) = \cx_R(\Omega(N)) $ (by Lemma~\ref{lem:cx-properties}(ii)),
	\begin{equation}\label{eqn1-thm:char-of-reg}
		\mbox{we may replace $ M $ and $ N $ by $ \Omega(M) $ and $ \Omega(N) $ respectively.}
	\end{equation}
	
	We proceed by induction on $ d := \dim(R) $. First assume that $ d = 0 $. In this case, we have $ \fm^3 = 0 $, and hence $ \fm^2 M = 0 = \fm^2 N $ by \eqref{eqn1-thm:char-of-reg} and Lemma~\ref{lemma: socle, syzygy}. If $ M $ or $ N $ is free, then $ \cx_R(k) = 0 $ (since $ \cx_R(k) = \cx_R(M) = \cx_R(N) $), which implies that $ \projdim_R(k) $ is finite, and hence $ R $ is regular. So we assume that $ M $ and $ N $ are not free. We show that this cannot be the case by getting a contradiction. Set $ \gamma(M) := \lambda(\fm M)/\mu(M) $. In view of \cite[Thm.~2.5(1) and (4)]{HSV04}, we obtain that $ \gamma(M) $ and $ \gamma(N) $ are positive integers satisfying $ \gamma(M)\gamma(N) = \type(R) = 1 $ (as $ R $ is Gorenstein), i.e., $ \gamma(M) = 1 = \gamma(N) $. It then follows from  \cite[Thm.~2.5(2)]{HSV04} that $ \beta_i^R(M) = \beta_0^R(M) $ $ (\neq 0) $ for all $ i \ge 1 $. Thus $ \cx_R(k) = \cx_R(M) = 1 $, which implies that $ R $ is a complete intersection ring of codimension $ c = 1 $; see Proposition~\ref{prop:facts-projcx-k}. Since $ \Tor_i^R(M,N) = 0 $ for all $ i \gg 0 $, and $ \cx_R(M) = \cx_R(N) = \cx_R(k) = c $, by virtue of \cite[Lem.~9.3.9]{Avr98}, it follows that $ c = 0 $, which is a contradiction.
	
	We now complete the inductive step as in the proof of Theorem~\ref{thm:Gor-rings-max-cx}. Assume that $ d \ge 1 $. Set $ \overline{(-)} := (-) \otimes_R R/(x_1)$, where $ \underline{x} = x_1,\ldots,x_d \in \fm \smallsetminus \fm^2 $. Replacing $ M $, $ N $ by $ \Omega(M) $, $ \Omega(N) $ respectively, we may assume that $ x_1 $ is regular on both $ M $ and $ N $. So, in view of Lemma~\ref{lem:Ext-Tor-mod-x}, $ \Tor_i^{\oR}(\oM,\oN) = 0 $ for all $ i \gg 0 $. By Lemmas~\ref{lem:cx-properties} and \ref{lem:cx-mod-x}, $ \cx_{\oR}(\oM) = \cx_R(M) = \cx_R(k) = \cx_{\oR}(k) $. Similarly, $ \cx_{\oR}(\oN) = \cx_{\oR}(k) $. Therefore, by induction hypothesis, $ \oR $ is regular, and hence $ R $ is regular.
\end{proof}

The following example shows that Theorem~\ref{thm:char-of-reg} does not necessarily hold true over arbitrary (i.e., non-Gorenstein) short local rings.

\begin{example}\label{exam1:thm:char-of-reg}
	Let $ R = k[x_1,\ldots,x_p,y_1,\ldots,y_q]/\left( (x_1,\ldots,x_p)^2 + (y_1,\ldots,y_q)^2 \right) $. Note that $ \fm^3 = 0 $, where $ \fm $ is the maximal ideal of $ R $. Set $ I := (x_1,\ldots,x_p) $, $ J := (y_1,\ldots,y_q) $, $ M := R/I $ and $ N := R/J $. It is well known that $ \Tor_1^R(R/I,R/J) \cong (I\cap J)/IJ $. Since $ I \cap J = IJ $, $ \Tor_1^R(M, N) = 0 $.	Note that $ M \cong (x_i) $ for every $ 1 \le i \le p $. Then $ \Omega_1^R(M) = I = (x_1) \oplus \cdots \oplus (x_p) \cong M^p $. By induction, $ \Omega_n^R(M) \cong M^{p^n} $ for every $ n \ge 1 $. Therefore
	\begin{align*}
		\Tor_n^R(M, N) & \cong \Tor_{n-1}^R\left(\Omega_1^R(M), N\right) \cong \cdots \cong \Tor_1^R\left(\Omega_{n-1}^R(M), N\right) \\
		& \cong \Tor_1^R\big(M^{p^{n-1}}, N\big) \cong \Tor_1^R(M, N)^{p^{n-1}} = 0
	\end{align*}
	for every $ n \ge 2 $. Following \cite[3.8(2)]{Les85}, we have $ \curv_R(M) = p $ $ ( = \curv_R(k) ) $ and $ \curv_R(N) = q $. We now assume that $ p \ge q \ge 2 $. It follows from Lemma~\ref{lem:cx-properties}(i)(b) that $ \cx_R(M) = \infty = \cx_R(N) $. Thus $ M $ and $ N $ have maximal complexity, and $ \Tor_n^R(M,N) = 0 $ for every $ n \ge 1 $, but $ R $ is not regular. Note that $ \Soc(R) = IJ $, which is not cyclic, hence $ R $ is not Gorenstein.
\end{example}
 
 Unlike Theorem~\ref{thm:Gor-rings-max-cx}, in Theorem~\ref{thm:char-of-reg}, the word `complexity' cannot be replaced by `curvature'.
 
\begin{example}\label{exam2:thm:char-of-reg}
	Let $ R = k[x,y]/(x^2,y^2) $. Clearly, $ R $ satisfies the hypotheses of Theorem~\ref{thm:char-of-reg}. Set $ M := R/(x) $ and $ N := R/(y) $. In view of Example~\ref{exam1:thm:char-of-reg}, $ M $ and $ N $ have maximal curvature, $ \Tor_i^R(M,N) = 0 $ for all $ i \ge 1 $, but $ R $ is not regular.
\end{example}

 As a consequence of Theorems~\ref{thm:Avr-max-cx} and \ref{thm:char-of-reg}, we obtain the following:

\begin{corollary}
	Let $ R $ be a Gorenstein short local ring.	Let $ M $ and $ N $ be nonzero homomorphic images of finite direct sums of syzygy modules of $ k $. $ ( $Possibly, $ M = N ) $. Then $ R $ is regular if and only if $ \Tor_i^R(M,N) = 0 $ for all $ i \gg 0 $.
\end{corollary}

\section*{acknowledgement}

The author is grateful to Anjan Gupta and Ryo Takahashi for making Lemma~\ref{lem:Gupta} and Example~\ref{exam:Takahashi} respectively. He is supported by IIT Hyderabad SEED Grant.

\end{document}